\documentclass[a4paper,reqno,12pt,oneside]{amsart}

\usepackage{amsfonts,amssymb,latexsym,xspace,epsfig,graphics,color}
\usepackage{amsmath,enumerate,stmaryrd,xy}
\usepackage{hyperref}
\usepackage{graphicx}
\usepackage{subfigure}

\usepackage{tikz}
\usetikzlibrary{arrows,snakes,backgrounds}

\def\D{\hbox{\sf D}}
\def\Stav{\hbox{\sf Stav}}

\def\Z{{\mathchoice
{\hbox{$\sf\textstyle         Z\kern-0.4em Z$}}
{\hbox{$\sf\textstyle         Z\kern-0.4em Z$}}
{\hbox{$\sf\scriptstyle       Z\kern-0.3em Z$}}
{\hbox{$\sf\scriptscriptstyle Z\kern-0.2em Z$}}
}}

\newtheorem{theorem}{Theorem}[section]

\title[An improved lower bound for the critical parameter of the Stavskaya's process]{An improved lower bound for the critical parameter\\ of the Stavskaya's process}

\author{Alex D. Ramos, Caliteia S. Sousa, Pablo M. Rodriguez and Paula Cadavid}

\address{
\newline
Alex D. Ramos; Pablo M. Rodriguez and Caliteia S. Sousa
\newline
Universidade Federal de Pernambuco, Department of Statistics, Recife, PE, 50740-540, Brazil.
\newline
e-mails: alex@de.ufpe.br; pablo@de.ufpe.br; caliteia@de.ufpe.br
\newline
\newline
Paula Cadavid
\newline
Universidade Federal do ABC, Santo Andr\'e, SP, 09210-580, Brazil.
\newline
e-mail: pacadavid@gmail.com
}


\begin{document}

\maketitle

\begin{abstract}
We consider the Stavskaya's process, which is a two-states Probabilistic Celular Automata defined on a one-dimensional lattice. The process is defined in such a way that the state of any vertex depends only on itself and on the state of its right-adjacent neighbor. This process was one of the first multicomponent systems with local interaction, for which has been proved rigorously the existence of a kind of phase transition. However, the exact localization of its critical value remains as an open problem. In this work we provide a new lower bound for the critical value. The last one was obtained by Andrei Toom, fifty years ago. 
\end{abstract}

\section{Introduction}

From the mid-twentieth century onwards, the development of a new part of the Theory of Stochastic Processes, called the Local Interaction Theory of Stochastic Processes, has begun to be developed, and now it is better known as the Theory of Interacting Particle Systems. The Stavskaya's process was one of the processes that contributed to this development. This process is a discrete-time version of the well-known contact process \cite{harris,Ligget:1985}, and may be described as a $\{0,1\}$-states Probabilistic Celular Automata (PCA) defined on a one-dimensional lattice. We assume that the state of any vertex depends only on itself and on the state of its right-adjacent neighbor. Moreover, each time of the process may be subdivided into a two-stages transition. In the first one each vertex of the lattice stay (or becomes) at state $1$ provided itself or its right-adjacent neighbor is at state $1$. On the other hand, in the second stage, each vertex at state $1$ turns $0$ with probability $\alpha$, independently of the other transitions. Thus defined, the constant $\alpha$ is the parameter of the model and it is associated with the randomness of the underlying stochastic process. It is not difficult to see that by considering the continuous-time version of this process we obtain the classical contact process. It is worth pointing out that although the contact process has been extensively studied in the literature, the Stavskaya's process has received (much) less attention and today it is an interesting resource of open problems. For some recent works dealing with existing open questions or generalizations for the Stavskaya's process we refer the reader to \cite{maes,Ponselet,Taggi}.

The Stavskaya's process is one of the first Interacting Particle Systems for which the existence of a phase transition has been rigorously proved \cite{ST2,ATC,Toom:2013,Toom:1990}.   More specifically, it was proved that there is $\alpha^{*}\in (0,1)$ such that for all $\alpha>\alpha^*$, the process is ergodic, i.e.,  the process started from any initial measure converge toward $\delta_0$. By another side, if $\alpha<\alpha^*,$ then the process started from the measure $\delta_1$ does not converge to $\delta_0.$ For this process, the exact value of $\alpha^*$ is not known, and only theoretical lower and upper bounds or estimates from computer simulations are available. A. Toom \cite{Toom:1990} proves that $\alpha^* \in (0.09, 0.323)$ and Mendon\c ca \cite{mendonca}, through computer simulations, estimates $\alpha^* \approx 0.29450(5).$

In this work, we revisited the method used in \cite{Toom:1990} to obtain the lower threshold to $\alpha^*,$ improving its estimation. We were able to offer now that $\alpha^* > 0.11$.

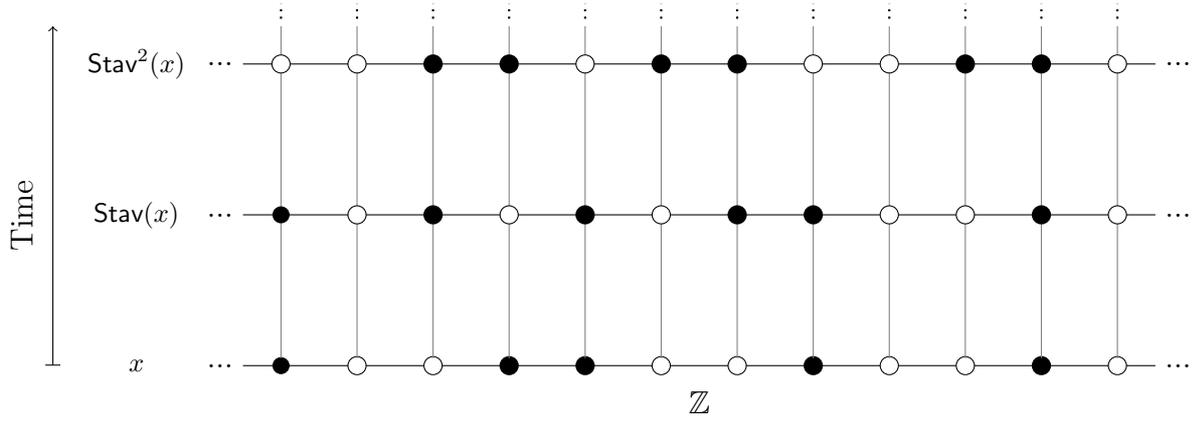
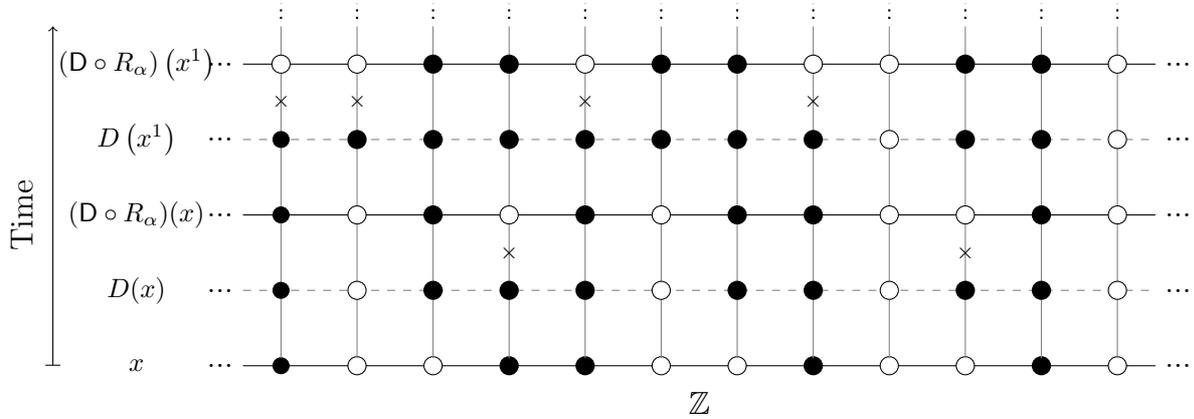
\begin{figure}[h]
\label{FIG:FP}
\begin{center}
\subfigure[][Realization of the Stavskaya's transformation with initial configuration $x$. The intermediate stages are represented in Figure 1(b).]{
\begin{tikzpicture}[scale=1]

\draw [|->] (-6,-2) -- (-6,2.5);
\draw [white] (10,-2) -- (10,2.5);
\draw (-6.7,0) node[below, rotate=90] {Time};
\draw (2.5,-2.2) node[below] {$\mathbb{Z}$};
\draw (-3.5,-2) -- (8.5,-2);
\draw (8.8,-2) node {...};
\draw (-3.8,-2) node {...};


\filldraw [black] (-3,-2) circle (3pt);
\draw [very thick] (-2,-2) circle (3pt);
\filldraw [white] (-2,-2) circle (3pt);
\draw [very thick] (-1,-2) circle (3pt);
\filldraw [white] (-1,-2) circle (3pt);
\draw [very thick] (0,-2) circle (3pt);
\filldraw [black] (0,-2) circle (3pt);
\draw [very thick] (1,-2) circle (3pt);
\filldraw [black] (1,-2) circle (3pt);
\draw [very thick] (2,-2) circle (3pt);
\filldraw [white] (2,-2) circle (3pt);
\draw [very thick] (3,-2) circle (3pt);
\filldraw [white] (3,-2) circle (3pt);
\draw [very thick] (4,-2) circle (3pt);
\filldraw [black] (4,-2) circle (3pt);
\draw [very thick] (5,-2) circle (3pt);
\filldraw [white] (5,-2) circle (3pt);
\draw [very thick] (6,-2) circle (3pt);
\filldraw [white] (6,-2) circle (3pt);
\draw [very thick] (7,-2) circle (3pt);
\filldraw [black] (7,-2) circle (3pt);
\draw [very thick] (8,-2) circle (3pt);
\filldraw [white] (8,-2) circle (3pt);

\draw [gray] (-3,-1.9) -- (-3,2.5);
\draw [gray] (-2,-1.9) -- (-2,2.5);
\draw [gray] (-1,-1.9) -- (-1,2.5);
\draw [gray] (0,-1.9) -- (0,2.5);
\draw [gray] (1,-1.9) -- (1,2.5);
\draw [gray] (2,-1.9) -- (2,2.5);
\draw [gray] (3,-1.9) -- (3,2.5);
\draw [gray] (4,-1.9) -- (4,2.5);
\draw [gray] (5,-1.9) -- (5,2.5);
\draw [gray] (6,-1.9) -- (6,2.5);
\draw [gray] (7,-1.9) -- (7,2.5);
\draw [gray] (8,-1.9) -- (8,2.5);

\draw [dotted, thick] (-3,2.6) -- (-3,2.8);
\draw [dotted, thick] (-2,2.6) -- (-2,2.8);
\draw [dotted, thick] (-1,2.6) -- (-1,2.8);
\draw [dotted, thick] (0,2.6) -- (0,2.8);
\draw [dotted, thick] (1,2.6) -- (1,2.8);
\draw [dotted, thick] (2,2.6) -- (2,2.8);
\draw [dotted, thick] (3,2.6) -- (3,2.8);
\draw [dotted, thick] (4,2.6) -- (4,2.8);
\draw [dotted, thick] (5,2.6) -- (5,2.8);
\draw [dotted, thick] (6,2.6) -- (6,2.8);
\draw [dotted, thick] (7,2.6) -- (7,2.8);
\draw [dotted, thick] (8,2.6) -- (8,2.8);



\draw (-3.5,0) -- (8.5,0);
\draw (8.8,0) node {...};
\draw (-3.8,0) node {...};

\filldraw [black] (-3,0) circle (3pt);
\draw [very thick] (-2,0) circle (3pt);
\filldraw [white] (-2,0) circle (3pt);
\draw [very thick] (-1,0) circle (3pt);
\filldraw [black] (-1,0) circle (3pt);
\draw [very thick] (0,0) circle (3pt);
\filldraw [white] (0,0) circle (3pt);
\draw [very thick] (1,0) circle (3pt);
\filldraw [black] (1,0) circle (3pt);
\draw [very thick] (2,0) circle (3pt);
\filldraw [white] (2,0) circle (3pt);
\draw [very thick] (3,0) circle (3pt);
\filldraw [black] (3,0) circle (3pt);
\draw [very thick] (4,0) circle (3pt);
\filldraw [black] (4,0) circle (3pt);
\draw [very thick] (5,0) circle (3pt);
\filldraw [white] (5,0) circle (3pt);
\draw [very thick] (6,0) circle (3pt);
\filldraw [white] (6,0) circle (3pt);
\draw [very thick] (7,0) circle (3pt);
\filldraw [black] (7,0) circle (3pt);
\draw [very thick] (8,0) circle (3pt);
\filldraw [white] (8,0) circle (3pt);




\draw (-3.5,2) -- (8.5,2);
\draw (8.8,2) node {...};
\draw (-3.8,2) node {...};

\draw [very thick] (-3,2) circle (3pt);
\filldraw [white] (-3,2) circle (3pt);
\draw [very thick] (-2,2) circle (3pt);
\filldraw [white] (-2,2) circle (3pt);
\draw [very thick] (-1,2) circle (3pt);
\filldraw [black] (-1,2) circle (3pt);
\draw [very thick] (0,2) circle (3pt);
\filldraw [black] (0,2) circle (3pt);
\draw [very thick] (1,2) circle (3pt);
\filldraw [white] (1,2) circle (3pt);
\draw [very thick] (2,2) circle (3pt);
\filldraw [black] (2,2) circle (3pt);
\draw [very thick] (3,2) circle (3pt);
\filldraw [black] (3,2) circle (3pt);
\draw [very thick] (4,2) circle (3pt);
\filldraw [white] (4,2) circle (3pt);
\draw [very thick] (5,2) circle (3pt);
\filldraw [white] (5,2) circle (3pt);
\draw [very thick] (6,2) circle (3pt);
\filldraw [black] (6,2) circle (3pt);
\draw [very thick] (7,2) circle (3pt);
\filldraw [black] (7,2) circle (3pt);
\draw [very thick] (8,2) circle (3pt);
\filldraw [white] (8,2) circle (3pt);

\draw (-4.9,-2) node[font=\footnotesize] {$x$};
\draw (-4.9,0) node[font=\footnotesize] {$\Stav(x)$};
\draw (-4.9,2) node[font=\footnotesize] {$\Stav^2(x)$};

\end{tikzpicture}}

\subfigure[][Realization of the Stavskaya's process stage by stage. The $\times$ marks are used to represent transitions from $1$ to $0$ coming from the $R_\alpha$ operator. The process start from a configuration $x$ and for simplicity we let $x^1:=(\D\circ R_{\alpha})(x)$.]{

\begin{tikzpicture}[scale=1]

\draw [|->] (-6,-2) -- (-6,2.5);
\draw [white] (10,-2) -- (10,2.5);
\draw (-6.7,0) node[below, rotate=90] {Time};
\draw (2.5,-2.2) node[below] {$\mathbb{Z}$};
\draw (-3.5,-2) -- (8.5,-2);
\draw (8.8,-2) node {...};
\draw (-3.8,-2) node {...};


\filldraw [black] (-3,-2) circle (3pt);
\draw [very thick] (-2,-2) circle (3pt);
\filldraw [white] (-2,-2) circle (3pt);
\draw [very thick] (-1,-2) circle (3pt);
\filldraw [white] (-1,-2) circle (3pt);
\draw [very thick] (0,-2) circle (3pt);
\filldraw [black] (0,-2) circle (3pt);
\draw [very thick] (1,-2) circle (3pt);
\filldraw [black] (1,-2) circle (3pt);
\draw [very thick] (2,-2) circle (3pt);
\filldraw [white] (2,-2) circle (3pt);
\draw [very thick] (3,-2) circle (3pt);
\filldraw [white] (3,-2) circle (3pt);
\draw [very thick] (4,-2) circle (3pt);
\filldraw [black] (4,-2) circle (3pt);
\draw [very thick] (5,-2) circle (3pt);
\filldraw [white] (5,-2) circle (3pt);
\draw [very thick] (6,-2) circle (3pt);
\filldraw [white] (6,-2) circle (3pt);
\draw [very thick] (7,-2) circle (3pt);
\filldraw [black] (7,-2) circle (3pt);
\draw [very thick] (8,-2) circle (3pt);
\filldraw [white] (8,-2) circle (3pt);

\draw [gray] (-3,-1.9) -- (-3,2.5);
\draw [gray] (-2,-1.9) -- (-2,2.5);
\draw [gray] (-1,-1.9) -- (-1,2.5);
\draw [gray] (0,-1.9) -- (0,2.5);
\draw [gray] (1,-1.9) -- (1,2.5);
\draw [gray] (2,-1.9) -- (2,2.5);
\draw [gray] (3,-1.9) -- (3,2.5);
\draw [gray] (4,-1.9) -- (4,2.5);
\draw [gray] (5,-1.9) -- (5,2.5);
\draw [gray] (6,-1.9) -- (6,2.5);
\draw [gray] (7,-1.9) -- (7,2.5);
\draw [gray] (8,-1.9) -- (8,2.5);

\draw [dotted, thick] (-3,2.6) -- (-3,2.8);
\draw [dotted, thick] (-2,2.6) -- (-2,2.8);
\draw [dotted, thick] (-1,2.6) -- (-1,2.8);
\draw [dotted, thick] (0,2.6) -- (0,2.8);
\draw [dotted, thick] (1,2.6) -- (1,2.8);
\draw [dotted, thick] (2,2.6) -- (2,2.8);
\draw [dotted, thick] (3,2.6) -- (3,2.8);
\draw [dotted, thick] (4,2.6) -- (4,2.8);
\draw [dotted, thick] (5,2.6) -- (5,2.8);
\draw [dotted, thick] (6,2.6) -- (6,2.8);
\draw [dotted, thick] (7,2.6) -- (7,2.8);
\draw [dotted, thick] (8,2.6) -- (8,2.8);


\draw [dashed, gray] (-3.5,-1) -- (8.5,-1);
\draw (8.8,-1) node {...};
\draw (-3.8,-1) node {...};

\filldraw [black] (-3,-1) circle (3pt);
\draw [very thick] (-2,-1) circle (3pt);
\filldraw [white] (-2,-1) circle (3pt);
\draw [very thick] (-1,-1) circle (3pt);
\filldraw [black] (-1,-1) circle (3pt);
\draw [very thick] (0,-1) circle (3pt);
\filldraw [black] (0,-1) circle (3pt);
\draw [very thick] (1,-1) circle (3pt);
\filldraw [black] (1,-1) circle (3pt);
\draw [very thick] (2,-1) circle (3pt);
\filldraw [white] (2,-1) circle (3pt);
\draw [very thick] (3,-1) circle (3pt);
\filldraw [black] (3,-1) circle (3pt);
\draw [very thick] (4,-1) circle (3pt);
\filldraw [black] (4,-1) circle (3pt);
\draw [very thick] (5,-1) circle (3pt);
\filldraw [white] (5,-1) circle (3pt);
\draw [very thick] (6,-1) circle (3pt);
\filldraw [black] (6,-1) circle (3pt);
\draw [very thick] (7,-1) circle (3pt);
\filldraw [black] (7,-1) circle (3pt);
\draw [very thick] (8,-1) circle (3pt);
\filldraw [white] (8,-1) circle (3pt);


\draw (-3.5,0) -- (8.5,0);
\draw (8.8,0) node {...};
\draw (-3.8,0) node {...};

\filldraw [black] (-3,0) circle (3pt);
\draw [very thick] (-2,0) circle (3pt);
\filldraw [white] (-2,0) circle (3pt);
\draw [very thick] (-1,0) circle (3pt);
\filldraw [black] (-1,0) circle (3pt);
\draw [very thick] (0,0) circle (3pt);
\filldraw [white] (0,0) circle (3pt);
\draw [very thick] (1,0) circle (3pt);
\filldraw [black] (1,0) circle (3pt);
\draw [very thick] (2,0) circle (3pt);
\filldraw [white] (2,0) circle (3pt);
\draw [very thick] (3,0) circle (3pt);
\filldraw [black] (3,0) circle (3pt);
\draw [very thick] (4,0) circle (3pt);
\filldraw [black] (4,0) circle (3pt);
\draw [very thick] (5,0) circle (3pt);
\filldraw [white] (5,0) circle (3pt);
\draw [very thick] (6,0) circle (3pt);
\filldraw [white] (6,0) circle (3pt);
\draw [very thick] (7,0) circle (3pt);
\filldraw [black] (7,0) circle (3pt);
\draw [very thick] (8,0) circle (3pt);
\filldraw [white] (8,0) circle (3pt);


\draw (0,-0.5) node[font=\tiny] {$\times$};
\draw (6,-0.5) node[font=\tiny] {$\times$};


\draw [dashed, gray] (-3.5,1) -- (8.5,1);
\draw (8.8,1) node {...};
\draw (-3.8,1) node {...};

\filldraw [black] (-3,1) circle (3pt);
\draw [very thick] (-2,1) circle (3pt);
\filldraw [black] (-2,1) circle (3pt);
\draw [very thick] (-1,1) circle (3pt);
\filldraw [black] (-1,1) circle (3pt);
\draw [very thick] (0,1) circle (3pt);
\filldraw [black] (0,1) circle (3pt);
\draw [very thick] (1,1) circle (3pt);
\filldraw [black] (1,1) circle (3pt);
\draw [very thick] (2,1) circle (3pt);
\filldraw [black] (2,1) circle (3pt);
\draw [very thick] (3,1) circle (3pt);
\filldraw [black] (3,1) circle (3pt);
\draw [very thick] (4,1) circle (3pt);
\filldraw [black] (4,1) circle (3pt);
\draw [very thick] (5,1) circle (3pt);
\filldraw [white] (5,1) circle (3pt);
\draw [very thick] (6,1) circle (3pt);
\filldraw [black] (6,1) circle (3pt);
\draw [very thick] (7,1) circle (3pt);
\filldraw [black] (7,1) circle (3pt);
\draw [very thick] (8,1) circle (3pt);
\filldraw [white] (8,1) circle (3pt);


\draw (-3.5,2) -- (8.5,2);
\draw (8.8,2) node {...};
\draw (-3.8,2) node {...};

\draw [very thick] (-3,2) circle (3pt);
\filldraw [white] (-3,2) circle (3pt);
\draw [very thick] (-2,2) circle (3pt);
\filldraw [white] (-2,2) circle (3pt);
\draw [very thick] (-1,2) circle (3pt);
\filldraw [black] (-1,2) circle (3pt);
\draw [very thick] (0,2) circle (3pt);
\filldraw [black] (0,2) circle (3pt);
\draw [very thick] (1,2) circle (3pt);
\filldraw [white] (1,2) circle (3pt);
\draw [very thick] (2,2) circle (3pt);
\filldraw [black] (2,2) circle (3pt);
\draw [very thick] (3,2) circle (3pt);
\filldraw [black] (3,2) circle (3pt);
\draw [very thick] (4,2) circle (3pt);
\filldraw [white] (4,2) circle (3pt);
\draw [very thick] (5,2) circle (3pt);
\filldraw [white] (5,2) circle (3pt);
\draw [very thick] (6,2) circle (3pt);
\filldraw [black] (6,2) circle (3pt);
\draw [very thick] (7,2) circle (3pt);
\filldraw [black] (7,2) circle (3pt);
\draw [very thick] (8,2) circle (3pt);
\filldraw [white] (8,2) circle (3pt);


\draw (-3,1.5) node[font=\tiny] {$\times$};
\draw (-2,1.5) node[font=\tiny] {$\times$};
\draw (1,1.5) node[font=\tiny] {$\times$};
\draw (4,1.5) node[font=\tiny] {$\times$};


\draw (-4.9,-2) node[font=\footnotesize] {$x$};
\draw (-4.9,-1) node[font=\footnotesize] {$D(x)$};
\draw (-4.9,0) node[font=\footnotesize] {$(\D\circ R_{\alpha})(x)$};
\draw (-4.9,1) node[font=\footnotesize] {$D\left(x^1\right)$};
\draw (-4.9,2) node[font=\footnotesize] {$(\D\circ R_{\alpha})\left(x^1\right)$};

\end{tikzpicture}}
\end{center}
\caption{Graphical space-time representation of the Stavskaya's process. Black and white particles represent vertices in state $1$ and $0$, respectively.}\label{fig:graphical}
\end{figure}

From now on, we shall describe the Stavskay's process in a formal way. Lets $\Z$ the set of integer numbers and we call $\{0,1\}^{\Z}$  the configuration space.  Every configuration $x$ is determined by its components $x_i\in\{1,0\}$, where $i\in\Z$. We shall consider a sequence of probabilistic measures enumerated by $t\in\{0,1,2,\ldots\}$, which we call {\it Stavskaya's process}. We assume that initially all the components are $1$ and then, at each step of the discrete time, two transformations occur. The first one is denoted by $\D$  and the second one is denoted by $R_{\alpha}$, where $\alpha\in[0,1].$ Thus, we define the Stavskaya's tranformation by $\Stav=\D\circ R_{\alpha}.$  
Speaking informally, when $\D$ is  applied to a configuration $x$ turns it into a configuration $y$ such that $y_i=max(x_{i},x_{i+1})$ for any $i\in \Z$, and  when $R_{\alpha}$ is applied it turns any $1$ at $0$ with probability $\alpha$, independently from what happens to other components. See Figure \ref{fig:graphical} for an illustration of a possible realization of this process. 

\bigskip

Now, we are able to declare our goal in this work:
\vspace{0.3cm}

\begin{theorem}\label{theor:main}
There is $\alpha^{*}>0.11$ such that if $\alpha<\alpha^*,$ then $\delta_1\Stav^t$ does not converge to $\delta_0$ when $t$ tends to infinity.
\end{theorem}

\section{Proof}

\subsection{A coupling }

\begin{figure}[!h]
\begin{center}
\subfigure[][Representation of $\Delta_{(0,8)}.$]{
\begin{tikzpicture}[scale=0.5]

\draw [->] (-8,0) -- (-7,1);
\draw [->] (-7,1) -- (-6,2);
\draw [->] (-6,2) -- (-5,3);
\draw [->] (-5,3) -- (-4,4);
\draw [->] (-4,4) -- (-3,5);
\draw [->] (-3,5) -- (-2,6);
\draw [->] (-2,6) -- (-1,7);
\draw [->] (-1,7) -- (0,8);
\draw [<-] (0,8) -- (1,7);
\draw [<-] (1,7) -- (2,6);
\draw [<-] (2,6) -- (3,5);
\draw [<-] (3,5) -- (4,4);
\draw [<-] (4,4) -- (5,3);
\draw [<-] (5,3) -- (6,2);
\draw [<-] (6,2) -- (7,1);
\draw [<-] (7,1) -- (8,0);

\draw [->] (-6,0) -- (-5,1);
\draw [->] (-5,1) -- (-4,2);
\draw [->] (-4,2) -- (-3,3);
\draw [->] (-3,3) -- (-2,4);
\draw [->] (-2,4) -- (-1,5);
\draw [->] (-1,5) -- (0,6);
\draw [<-] (0,6) -- (1,5);
\draw [<-] (1,5) -- (2,4);
\draw [<-] (2,4) -- (3,3);
\draw [<-] (3,3) -- (4,2);
\draw [<-] (4,2) -- (5,1);
\draw [<-] (5,1) -- (6,0);

\draw [->] (-4,0) -- (-3,1);
\draw [->] (-3,1) -- (-2,2);
\draw [->] (-2,2) -- (-1,3);
\draw [->] (-1,3) -- (0,4);
\draw [<-] (0,4) -- (1,3);
\draw [<-] (1,3) -- (2,2);
\draw [<-] (2,2) -- (3,1);
\draw [<-] (3,1) -- (4,0);

\draw [->] (-2,0) -- (-1,1);
\draw [->] (-1,1) -- (0,2);
\draw [<-] (0,2) -- (1,1);
\draw [<-] (1,1) -- (2,0);

\draw [->] (-6,0) -- (-7,1);
\draw [->] (-4,0) -- (-5,1);
\draw [->] (-5,1) -- (-6,2);
\draw [->] (-2,0) -- (-3,1);
\draw [->] (-3,1) -- (-4,2);
\draw [->] (-4,2) -- (-5,3);
\draw [->] (0,0) -- (-1,1);
\draw [->] (-1,1) -- (-2,2);
\draw [->] (-2,2) -- (-3,3);
\draw [->] (-3,3) -- (-4,4);
\draw [->] (0,2) -- (-1,3);
\draw [->] (-1,3) -- (-2,4);
\draw [->] (-2,4) -- (-3,5);
\draw [->] (0,4) -- (-1,5);
\draw [->] (-1,5) -- (-2,6);
\draw [->] (0,6) -- (-1,7);
\draw [->] (0,6) -- (1,7);
\draw [->] (0,4) -- (1,5);
\draw [->] (1,5) -- (2,6);
\draw [->] (0,2) -- (1,3);
\draw [->] (1,3) -- (2,4);
\draw [->] (2,4) -- (3,5);
\draw [->] (0,0) -- (1,1);
\draw [->] (1,1) -- (2,2);
\draw [->] (2,2) -- (3,3);
\draw [->] (3,3) -- (4,4);
\draw [->] (2,0) -- (3,1);
\draw [->] (3,1) -- (4,2);
\draw [->] (4,2) -- (5,3);
\draw [->] (4,0) -- (5,1);
\draw [->] (5,1) -- (6,2);
\draw [->] (6,0) -- (7,1);
\end{tikzpicture}}\qquad\qquad\qquad \subfigure[][The orientated bond percolation graph, for $\Delta_{(0,8)}$, after the replacement of vertices by vertical edges.]{
\begin{tikzpicture}[scale=0.35]

\draw [->] (-8,-1) -- (-8,0);
\draw [->] (-6,-1) -- (-6,0);
\draw [->] (-4,-1) -- (-4,0);
\draw [->] (-2,-1) -- (-2,0);
\draw [->] (0,-1) -- (0,0);
\draw [->] (2,-1) -- (2,0);
\draw [->] (4,-1) -- (4,0);
\draw [->] (6,-1) -- (6,0);
\draw [->] (8,-1) -- (8,0);

\draw (0,-3) -- (8,-1);
\draw (0,-3) -- (6,-1);
\draw (0,-3) -- (4,-1);
\draw (0,-3) -- (2,-1);
\draw (0,-3) -- (0,-1);
\draw (0,-3) -- (-2,-1);
\draw (0,-3) -- (-4,-1);
\draw (0,-3) -- (-6,-1);
\draw (0,-3) -- (-8,-1);

\draw [very thick] (0,-3) circle (3pt);
\filldraw [black] (0,-3) circle (3pt);
\draw (0,-3.8) node[font=\small] {$F$};

\draw [->] (-7,1) -- (-7,2);
\draw [->] (-5,1) -- (-5,2);
\draw [->] (-3,1) -- (-3,2);
\draw [->] (-1,1) -- (-1,2);
\draw [->] (1,1) -- (1,2);
\draw [->] (3,1) -- (3,2);
\draw [->] (5,1) -- (5,2);
\draw [->] (7,1) -- (7,2);

\draw [->] (-6,3) -- (-6,4);
\draw [->] (-4,3) -- (-4,4);
\draw [->] (-2,3) -- (-2,4);
\draw [->] (0,3) -- (0,4);
\draw [->] (2,3) -- (2,4);
\draw [->] (4,3) -- (4,4);
\draw [->] (6,3) -- (6,4);

\draw [->] (-5,5) -- (-5,6);
\draw [->] (-3,5) -- (-3,6);
\draw [->] (-1,5) -- (-1,6);
\draw [->] (1,5) -- (1,6);
\draw [->] (3,5) -- (3,6);
\draw [->] (5,5) -- (5,6);

\draw [->] (-4,7) -- (-4,8);
\draw [->] (-2,7) -- (-2,8);
\draw [->] (0,7) -- (0,8);
\draw [->] (2,7) -- (2,8);
\draw [->] (4,7) -- (4,8);

\draw [->] (-3,9) -- (-3,10);
\draw [->] (-1,9) -- (-1,10);
\draw [->] (1,9) -- (1,10);
\draw [->] (3,9) -- (3,10);

\draw [->] (-2,11) -- (-2,12);
\draw [->] (0,11) -- (0,12);
\draw [->] (2,11) -- (2,12);

\draw [->] (-1,13) -- (-1,14);
\draw [->] (1,13) -- (1,14);

\draw [->] (0,15) -- (0,16);

\draw [->] (-8,0) -- (-7,1);
\draw [->] (-7,2) -- (-6,3);
\draw [->] (-6,4) -- (-5,5);
\draw [->] (-5,6) -- (-4,7);
\draw [->] (-4,8) -- (-3,9);
\draw [->] (-3,10) -- (-2,11);
\draw [->] (-2,12) -- (-1,13);
\draw [->] (-1,14) -- (0,15);
\draw [<-] (0,15) -- (1,14);
\draw [<-] (1,13) -- (2,12);
\draw [<-] (2,11) -- (3,10);
\draw [<-] (3,9) -- (4,8);
\draw [<-] (4,7) -- (5,6);
\draw [<-] (5,5) -- (6,4);
\draw [<-] (6,3) -- (7,2);
\draw [<-] (7,1) -- (8,0);

\draw [->] (-6,0) -- (-5,1);
\draw [->] (-5,2) -- (-4,3);
\draw [->] (-4,4) -- (-3,5);
\draw [->] (-3,6) -- (-2,7);
\draw [->] (-2,8) -- (-1,9);
\draw [->] (-1,10) -- (0,11);
\draw [<-] (0,11) -- (1,10);
\draw [<-] (1,9) -- (2,8);
\draw [<-] (2,7) -- (3,6);
\draw [<-] (3,5) -- (4,4);
\draw [<-] (4,3) -- (5,2);
\draw [<-] (5,1) -- (6,0);

\draw [->] (-4,0) -- (-3,1);
\draw [->] (-3,2) -- (-2,3);
\draw [->] (-2,4) -- (-1,5);
\draw [->] (-1,6) -- (0,7);
\draw [<-] (0,7) -- (1,6);
\draw [<-] (1,5) -- (2,4);
\draw [<-] (2,3) -- (3,2);
\draw [<-] (3,1) -- (4,0);

\draw [->] (-2,0) -- (-1,1);
\draw [->] (-1,2) -- (0,3);
\draw [<-] (0,3) -- (1,2);
\draw [<-] (1,1) -- (2,0);

\draw [->] (-6,0) -- (-7,1);
\draw [->] (-4,0) -- (-5,1);
\draw [->] (-5,2) -- (-6,3);
\draw [->] (-2,0) -- (-3,1);
\draw [->] (-3,2) -- (-4,3);
\draw [->] (-4,4) -- (-5,5);
\draw [->] (0,0) -- (-1,1);
\draw [->] (-1,2) -- (-2,3);
\draw [->] (-2,4) -- (-3,5);
\draw [->] (-3,6) -- (-4,7);
\draw [->] (0,4) -- (-1,5);
\draw [->] (-1,6) -- (-2,7);
\draw [->] (-2,8) -- (-3,9);
\draw [->] (0,12) -- (-1,13);
\draw [->] (-1,10) -- (-2,11);
\draw [->] (0,8) -- (-1,9);
\draw [->] (0,8) -- (1,9);
\draw [->] (0,12) -- (1,13);
\draw [->] (1,10) -- (2,11);
\draw [->] (0,4) -- (1,5);
\draw [->] (1,6) -- (2,7);
\draw [->] (2,8) -- (3,9);
\draw [->] (0,0) -- (1,1);
\draw [->] (1,2) -- (2,3);
\draw [->] (2,4) -- (3,5);
\draw [->] (3,6) -- (4,7);
\draw [->] (2,0) -- (3,1);
\draw [->] (3,2) -- (4,3);
\draw [->] (4,4) -- (5,5);
\draw [->] (4,0) -- (5,1);
\draw [->] (5,2) -- (6,3);
\draw [->] (6,0) -- (7,1);
\end{tikzpicture}}
	\caption{{\it Illustration of the oriented percolation graph.}}\label{fig1}
\end{center}

\end{figure}
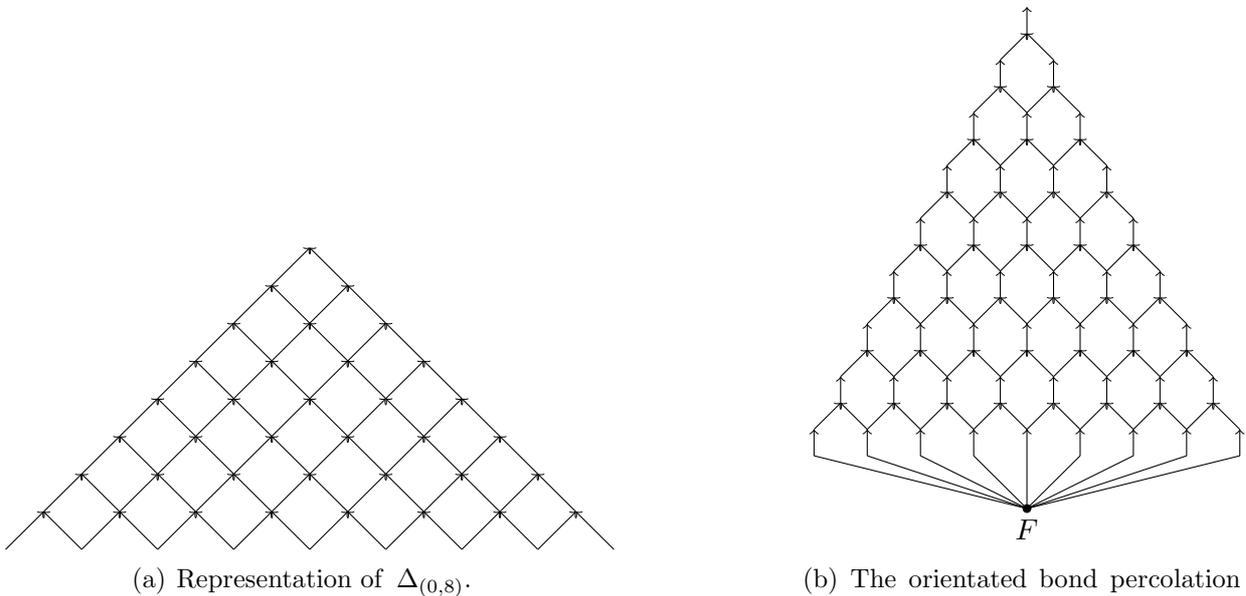

It is well-known that the Stavskaya's process can be represented using oriented percolation. In our work, by technical question, we shall define this oriented graph slightly different from that ordinary used. The  vertices of our planar graph are given by  $(i,t)$ where $i\in\Z$ and $t\in\Z_+.$ From every vertex $(i,t+1)$ oriented bonds come from the vertices $(i,t)$ and $(i+1,t).$ Every bond from $t$ to $t+1$ is open and closed in the opposite direction, that is from $t+1$ to $t$. We assume that all the vertices $(i,0)$ are open and the others vertices are closed with probability $\alpha$ and open with probability $1-\alpha$ independently of what occur on the other vertices. Notice that a vertex $(i,t)$ has a particle only if there is an open path (i.e. formed by open bonds and vertices) on the set
\[
\Delta_{(i,t)}=\{(j,s): 0\leq s \leq t \text{ and }i \le j \le i+t-s\},
\]  
connecting some vertex $(j,0)$ for $i\leq j \leq i+t$ with the vertex $(i,t)$. Note that for a fixed vertex $(i,t),$ the $\Delta_{(i,t)}$ has $(t+2)(t+1)/2$ vertices. To get a better graphical representation with a triangle, the vertical bonds are inclined to the right-side. See Figure \ref{fig1}(a), for an illustration of a fragment of the percolation graph for the triangle $\Delta_{(0,8)}.$  

In \cite{Toom:1990} it has been proved that in the Stavskaya's process, under the assumption adopted here, there is a particle at a component on the position $i$ at time $t$ if, and only if, there is a open path from some vertex $(j,0)$ to the vertex $(i,t)$ in the oriented percolation. This result is proved by mean of the coupling between the Stavskaya's process and the oriented percolation model, namely: one can consider the states of the initial vertices, all open, with the states of the  configuration where all the particles are on the state $1$, whose the initial measure is concentrated. The other vertices assume the state, open or closed, according with the action of the operator $R_{\alpha}.$ By end, the inclined bonds are associated with the action of the operator $\D.$    

We shall change our oriented site-bond percolation by a oriented bond percolation, for this we replace all the vertices by a vertical bond orientated upper ward, which will be open or closed following the same previous assumption from its corresponding vertex. 
Each bond corresponding to an initial vertex, $(i,0),$ is analog a one font. Thus, we shall have a great quantity of fonts. To avoid it, we  establish only one font, which we will denote by $F$ and will be connected with open bonds at both directions to all initial bonds. On the Figure \ref{fig1}(b), we exhibits the fragment of the percolation graph, for $\Delta_{(0,8)},$ after this replacement and insertion of this font. 

\begin{figure}
\begin{center}
\subfigure[][The orientated percolation graph of $\Delta_{(0,8)}$ and its respective dual graph.]{
\begin{tikzpicture}[scale=0.4]

\draw [->] (-8,-1) -- (-8,0);
\draw [->] (-6,-1) -- (-6,0);
\draw [->] (-4,-1) -- (-4,0);
\draw [->] (-2,-1) -- (-2,0);
\draw [->] (0,-1) -- (0,0);
\draw [->] (2,-1) -- (2,0);
\draw [->] (4,-1) -- (4,0);
\draw [->] (6,-1) -- (6,0);
\draw [->] (8,-1) -- (8,0);

\draw (0,-3) -- (8,-1);
\draw (0,-3) -- (6,-1);
\draw (0,-3) -- (4,-1);
\draw (0,-3) -- (2,-1);
\draw (0,-3) -- (0,-1);
\draw (0,-3) -- (-2,-1);
\draw (0,-3) -- (-4,-1);
\draw (0,-3) -- (-6,-1);
\draw (0,-3) -- (-8,-1);

\draw [dashed, gray] (-9,-0.5) -- (9,-0.5);
\draw [dashed, gray] (-8,1.5) -- (8,1.5);
\draw [dashed, gray] (-7,3.5) -- (7,3.5);
\draw [dashed, gray] (-6,5.5) -- (6,5.5);
\draw [dashed, gray] (-5,7.5) -- (5,7.5);
\draw [dashed, gray] (-4,9.5) -- (4,9.5);
\draw [dashed, gray] (-3,11.5) -- (3,11.5);
\draw [dashed, gray] (-2,13.5) -- (2,13.5);
\draw [dashed, gray] (-1,15.5) -- (1,15.5);

\draw [dashed, gray] (-9,-0.5) -- (0,17.5) -- (9,-0.5);
\draw [dashed, gray] (-7,-0.5) -- (-8,1.5);
\draw [dashed, gray] (-5,-0.5) -- (-7,3.5);
\draw [dashed, gray] (-3,-0.5) -- (-6,5.5);
\draw [dashed, gray] (-1,-0.5) -- (-5,7.5);
\draw [dashed, gray] (1,-0.5) -- (-4,9.5);
\draw [dashed, gray] (3,-0.5) -- (-3,11.5);
\draw [dashed, gray] (5,-0.5) -- (-2,13.5);
\draw [dashed, gray] (7,-0.5) -- (-1,15.5);

\draw [dashed, gray] (7,-0.5) -- (8,1.5);
\draw [dashed, gray] (5,-0.5) -- (7,3.5);
\draw [dashed, gray] (3,-0.5) -- (6,5.5);
\draw [dashed, gray] (1,-0.5) -- (5,7.5);
\draw [dashed, gray] (-1,-0.5) -- (4,9.5);
\draw [dashed, gray] (-3,-0.5) -- (3,11.5);
\draw [dashed, gray] (-5,-0.5) -- (2,13.5);
\draw [dashed, gray] (-7,-0.5) -- (1,15.5);

\draw [very thick] (0,-3) circle (3pt);
\filldraw [black] (0,-3) circle (3pt);

\draw [->] (-7,1) -- (-7,2);
\draw [->] (-5,1) -- (-5,2);
\draw [->] (-3,1) -- (-3,2);
\draw [->] (-1,1) -- (-1,2);
\draw [->] (1,1) -- (1,2);
\draw [->] (3,1) -- (3,2);
\draw [->] (5,1) -- (5,2);
\draw [->] (7,1) -- (7,2);

\draw [->] (-6,3) -- (-6,4);
\draw [->] (-4,3) -- (-4,4);
\draw [->] (-2,3) -- (-2,4);
\draw [->] (0,3) -- (0,4);
\draw [->] (2,3) -- (2,4);
\draw [->] (4,3) -- (4,4);
\draw [->] (6,3) -- (6,4);

\draw [->] (-5,5) -- (-5,6);
\draw [->] (-3,5) -- (-3,6);
\draw [->] (-1,5) -- (-1,6);
\draw [->] (1,5) -- (1,6);
\draw [->] (3,5) -- (3,6);
\draw [->] (5,5) -- (5,6);

\draw [->] (-4,7) -- (-4,8);
\draw [->] (-2,7) -- (-2,8);
\draw [->] (0,7) -- (0,8);
\draw [->] (2,7) -- (2,8);
\draw [->] (4,7) -- (4,8);

\draw [->] (-3,9) -- (-3,10);
\draw [->] (-1,9) -- (-1,10);
\draw [->] (1,9) -- (1,10);
\draw [->] (3,9) -- (3,10);

\draw [->] (-2,11) -- (-2,12);
\draw [->] (0,11) -- (0,12);
\draw [->] (2,11) -- (2,12);

\draw [->] (-1,13) -- (-1,14);
\draw [->] (1,13) -- (1,14);

\draw [->] (0,15) -- (0,16);

\draw [->] (-8,0) -- (-7,1);
\draw [->] (-7,2) -- (-6,3);
\draw [->] (-6,4) -- (-5,5);
\draw [->] (-5,6) -- (-4,7);
\draw [->] (-4,8) -- (-3,9);
\draw [->] (-3,10) -- (-2,11);
\draw [->] (-2,12) -- (-1,13);
\draw [->] (-1,14) -- (0,15);
\draw [<-] (0,15) -- (1,14);
\draw [<-] (1,13) -- (2,12);
\draw [<-] (2,11) -- (3,10);
\draw [<-] (3,9) -- (4,8);
\draw [<-] (4,7) -- (5,6);
\draw [<-] (5,5) -- (6,4);
\draw [<-] (6,3) -- (7,2);
\draw [<-] (7,1) -- (8,0);

\draw [->] (-6,0) -- (-5,1);
\draw [->] (-5,2) -- (-4,3);
\draw [->] (-4,4) -- (-3,5);
\draw [->] (-3,6) -- (-2,7);
\draw [->] (-2,8) -- (-1,9);
\draw [->] (-1,10) -- (0,11);
\draw [<-] (0,11) -- (1,10);
\draw [<-] (1,9) -- (2,8);
\draw [<-] (2,7) -- (3,6);
\draw [<-] (3,5) -- (4,4);
\draw [<-] (4,3) -- (5,2);
\draw [<-] (5,1) -- (6,0);

\draw [->] (-4,0) -- (-3,1);
\draw [->] (-3,2) -- (-2,3);
\draw [->] (-2,4) -- (-1,5);
\draw [->] (-1,6) -- (0,7);
\draw [<-] (0,7) -- (1,6);
\draw [<-] (1,5) -- (2,4);
\draw [<-] (2,3) -- (3,2);
\draw [<-] (3,1) -- (4,0);

\draw [->] (-2,0) -- (-1,1);
\draw [->] (-1,2) -- (0,3);
\draw [<-] (0,3) -- (1,2);
\draw [<-] (1,1) -- (2,0);

\draw [->] (-6,0) -- (-7,1);
\draw [->] (-4,0) -- (-5,1);
\draw [->] (-5,2) -- (-6,3);
\draw [->] (-2,0) -- (-3,1);
\draw [->] (-3,2) -- (-4,3);
\draw [->] (-4,4) -- (-5,5);
\draw [->] (0,0) -- (-1,1);
\draw [->] (-1,2) -- (-2,3);
\draw [->] (-2,4) -- (-3,5);
\draw [->] (-3,6) -- (-4,7);
\draw [->] (0,4) -- (-1,5);
\draw [->] (-1,6) -- (-2,7);
\draw [->] (-2,8) -- (-3,9);
\draw [->] (0,12) -- (-1,13);
\draw [->] (-1,10) -- (-2,11);
\draw [->] (0,8) -- (-1,9);
\draw [->] (0,8) -- (1,9);
\draw [->] (0,12) -- (1,13);
\draw [->] (1,10) -- (2,11);
\draw [->] (0,4) -- (1,5);
\draw [->] (1,6) -- (2,7);
\draw [->] (2,8) -- (3,9);
\draw [->] (0,0) -- (1,1);
\draw [->] (1,2) -- (2,3);
\draw [->] (2,4) -- (3,5);
\draw [->] (3,6) -- (4,7);
\draw [->] (2,0) -- (3,1);
\draw [->] (3,2) -- (4,3);
\draw [->] (4,4) -- (5,5);
\draw [->] (4,0) -- (5,1);
\draw [->] (5,2) -- (6,3);
\draw [->] (6,0) -- (7,1);
\end{tikzpicture}}\qquad\qquad\qquad\subfigure[][The subgraph of $\Delta_{(0,8)}$ which is a trapezoid.]{
\begin{tikzpicture}[scale=0.4]

\draw [->] (-8,-1) -- (-8,0);
\draw [->] (-6,-1) -- (-6,0);
\draw [->] (-4,-1) -- (-4,0);
\draw [->] (-2,-1) -- (-2,0);
\draw [->] (0,-1) -- (0,0);
\draw [->] (2,-1) -- (2,0);
\draw [->] (4,-1) -- (4,0);
\draw [->] (6,-1) -- (6,0);
\draw [->] (8,-1) -- (8,0);

\draw (0,-3) -- (8,-1);
\draw (0,-3) -- (6,-1);
\draw (0,-3) -- (4,-1);
\draw (0,-3) -- (2,-1);
\draw (0,-3) -- (0,-1);
\draw (0,-3) -- (-2,-1);
\draw (0,-3) -- (-4,-1);
\draw (0,-3) -- (-6,-1);
\draw (0,-3) -- (-8,-1);

\draw (0,12) -- (1,10);
\draw (0,12) -- (3,10);
\draw (0,12) -- (-1,10);
\draw (0,12) -- (-3,10);

\draw [dashed, thick] (-9,-0.5) -- (9,-0.5);
\draw [dashed, gray] (-8,1.5) -- (8,1.5);
\draw [dashed, gray] (-7,3.5) -- (7,3.5);
\draw [dashed, gray] (-6,5.5) -- (6,5.5);
\draw [dashed, gray] (-5,7.5) -- (5,7.5);
\draw [dashed, thick] (-4,9.5) -- (4,9.5);

\draw [dashed, thick] (-9,-0.5) -- (-4,9.5);
\draw [dashed, thick] (9,-0.5) -- (4,9.5);
\draw [dashed, gray] (-7,-0.5) -- (-8,1.5);
\draw [dashed, gray] (-5,-0.5) -- (-7,3.5);
\draw [dashed, gray] (-3,-0.5) -- (-6,5.5);
\draw [dashed, gray] (-1,-0.5) -- (-5,7.5);
\draw [dashed, gray] (1,-0.5) -- (-4,9.5);
\draw [dashed, gray] (3,-0.5) -- (-2,9.5);
\draw [dashed, gray] (5,-0.5) -- (0,9.5);
\draw [dashed, gray] (7,-0.5) -- (2,9.5);

\draw [dashed, gray] (7,-0.5) -- (8,1.5);
\draw [dashed, gray] (5,-0.5) -- (7,3.5);
\draw [dashed, gray] (3,-0.5) -- (6,5.5);
\draw [dashed, gray] (1,-0.5) -- (5,7.5);
\draw [dashed, gray] (-1,-0.5) -- (4,9.5);
\draw [dashed, gray] (-3,-0.5) -- (2,9.5);
\draw [dashed, gray] (-5,-0.5) -- (0,9.5);
\draw [dashed, gray] (-7,-0.5) -- (-2,9.5);

\draw [very thick] (0,-3) circle (3pt);
\filldraw [black] (0,-3) circle (3pt);

\draw [very thick] (0,12) circle (3pt);
\filldraw [black] (0,12) circle (3pt);

\draw [->] (-7,1) -- (-7,2);
\draw [->] (-5,1) -- (-5,2);
\draw [->] (-3,1) -- (-3,2);
\draw [->] (-1,1) -- (-1,2);
\draw [->] (1,1) -- (1,2);
\draw [->] (3,1) -- (3,2);
\draw [->] (5,1) -- (5,2);
\draw [->] (7,1) -- (7,2);

\draw [->] (-6,3) -- (-6,4);
\draw [->] (-4,3) -- (-4,4);
\draw [->] (-2,3) -- (-2,4);
\draw [->] (0,3) -- (0,4);
\draw [->] (2,3) -- (2,4);
\draw [->] (4,3) -- (4,4);
\draw [->] (6,3) -- (6,4);

\draw [->] (-5,5) -- (-5,6);
\draw [->] (-3,5) -- (-3,6);
\draw [->] (-1,5) -- (-1,6);
\draw [->] (1,5) -- (1,6);
\draw [->] (3,5) -- (3,6);
\draw [->] (5,5) -- (5,6);

\draw [->] (-4,7) -- (-4,8);
\draw [->] (-2,7) -- (-2,8);
\draw [->] (0,7) -- (0,8);
\draw [->] (2,7) -- (2,8);
\draw [->] (4,7) -- (4,8);

\draw [->] (-3,9) -- (-3,10);
\draw [->] (-1,9) -- (-1,10);
\draw [->] (1,9) -- (1,10);
\draw [->] (3,9) -- (3,10);

\draw [->] (-8,0) -- (-7,1);
\draw [->] (-7,2) -- (-6,3);
\draw [->] (-6,4) -- (-5,5);
\draw [->] (-5,6) -- (-4,7);
\draw [->] (-4,8) -- (-3,9);
\draw [<-] (3,9) -- (4,8);
\draw [<-] (4,7) -- (5,6);
\draw [<-] (5,5) -- (6,4);
\draw [<-] (6,3) -- (7,2);
\draw [<-] (7,1) -- (8,0);

\draw [->] (-6,0) -- (-5,1);
\draw [->] (-5,2) -- (-4,3);
\draw [->] (-4,4) -- (-3,5);
\draw [->] (-3,6) -- (-2,7);
\draw [->] (-2,8) -- (-1,9);
\draw [<-] (1,9) -- (2,8);
\draw [<-] (2,7) -- (3,6);
\draw [<-] (3,5) -- (4,4);
\draw [<-] (4,3) -- (5,2);
\draw [<-] (5,1) -- (6,0);

\draw [->] (-4,0) -- (-3,1);
\draw [->] (-3,2) -- (-2,3);
\draw [->] (-2,4) -- (-1,5);
\draw [->] (-1,6) -- (0,7);
\draw [<-] (0,7) -- (1,6);
\draw [<-] (1,5) -- (2,4);
\draw [<-] (2,3) -- (3,2);
\draw [<-] (3,1) -- (4,0);

\draw [->] (-2,0) -- (-1,1);
\draw [->] (-1,2) -- (0,3);
\draw [<-] (0,3) -- (1,2);
\draw [<-] (1,1) -- (2,0);

\draw [->] (-6,0) -- (-7,1);
\draw [->] (-4,0) -- (-5,1);
\draw [->] (-5,2) -- (-6,3);
\draw [->] (-2,0) -- (-3,1);
\draw [->] (-3,2) -- (-4,3);
\draw [->] (-4,4) -- (-5,5);
\draw [->] (0,0) -- (-1,1);
\draw [->] (-1,2) -- (-2,3);
\draw [->] (-2,4) -- (-3,5);
\draw [->] (-3,6) -- (-4,7);
\draw [->] (0,4) -- (-1,5);
\draw [->] (-1,6) -- (-2,7);
\draw [->] (-2,8) -- (-3,9);
\draw [->] (0,8) -- (-1,9);
\draw [->] (0,8) -- (1,9);
\draw [->] (0,4) -- (1,5);
\draw [->] (1,6) -- (2,7);
\draw [->] (2,8) -- (3,9);
\draw [->] (0,0) -- (1,1);
\draw [->] (1,2) -- (2,3);
\draw [->] (2,4) -- (3,5);
\draw [->] (3,6) -- (4,7);
\draw [->] (2,0) -- (3,1);
\draw [->] (3,2) -- (4,3);
\draw [->] (4,4) -- (5,5);
\draw [->] (4,0) -- (5,1);
\draw [->] (5,2) -- (6,3);
\draw [->] (6,0) -- (7,1);
\end{tikzpicture}}
\caption{Representation for the oriented percolation graph associated to $\Delta_{(0,8)}$ and its respective dual graph.}\label{fig2}
\end{center}
\end{figure}
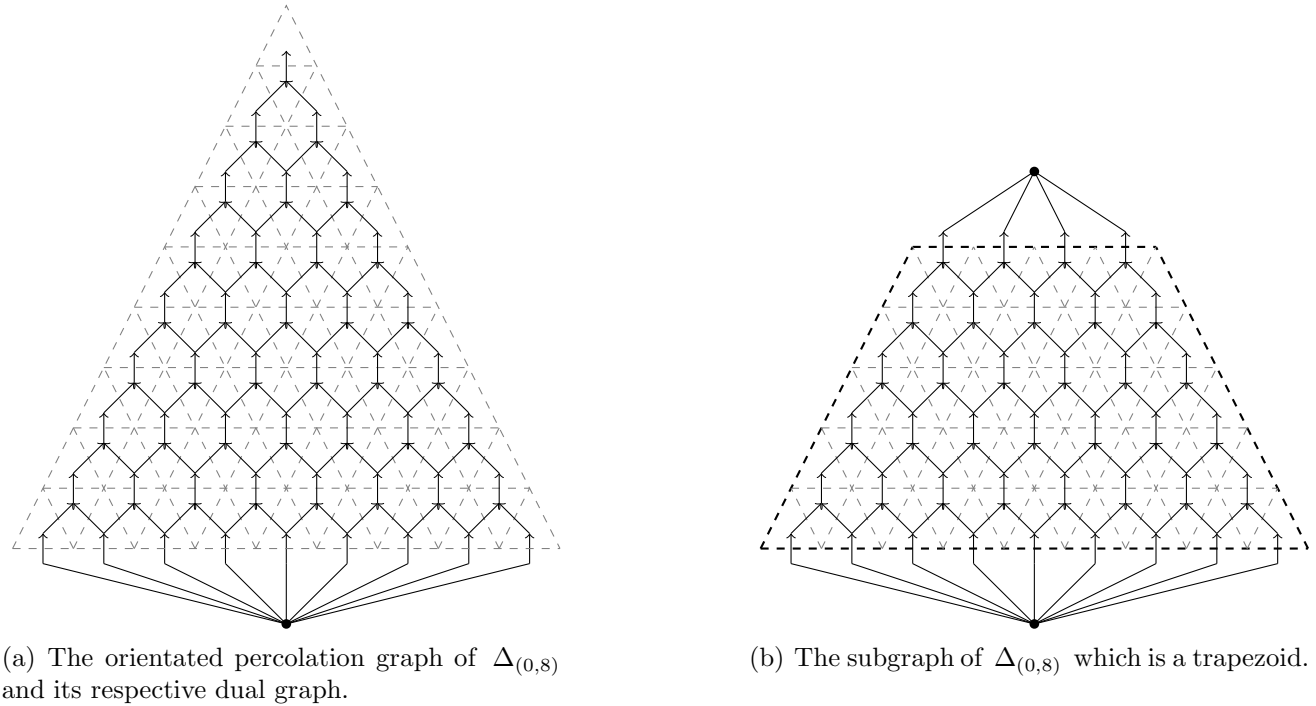

As usual, we can define the respective dual graph, see Figure \ref{fig2}. In order to do it, we consider our oriented graph percolation, and we get the dual graph by assuming that the directed bonds $\swarrow$ and $\nwarrow$ are always open and closed in the opposite directions, while the directed bonds $\longrightarrow$ are open with probability $\alpha$ and closed in the opposite direction. 

It is a well-known fact that there is no percolation in the original graph if there is a open contour going in the counterclockwise direction at the dual graph, surround the vertical bond corresponding to the ``peak" vertex of $\Delta_{(i,t)},$ the $(i,t).$  Through the coupling between the Stavskaya's process and the oriented percolation model, if the probability that there is such contour is less than $1$, then the $\delta_1\Stav^t$ does not converge to $\delta_0$ when $t\to\infty.$ It is what we will do to prove Theorem \ref{theor:main}.

\subsection{The recurrent method}

Lets $C_k$ the number of  contours with $k$ horizontal bonds on the dual graph, which starts from the left border  to the right border of the trapezoid (see Figure \ref{fig2}(b)). It was proved in \cite{Toom:1990} that
\begin{equation}
\delta_1\Stav^t(x_1=0,\ldots,x_m=0)\le\sum_{k=1}^mC_k\alpha^k.
\label{eq:limitesup}
\end{equation}
When the quantity in the left side of (\ref{eq:limitesup}) is equal to one, we have that in the corresponding oriented graph percolation there is a barrier, which will not permit the percolation. So, to prove Theorem \ref{theor:main} it is sufficient  to verify when the right-side of (\ref{eq:limitesup}) is less than one, and it is what we will do.  

Let us consider a coordinate system in the Figure \ref{fig2} right-side, whose the origin of the system is the left-upper border of the trapezoid. Given a vertex in the graph, the {\it Shift } for to do the contour is a two-dimensional vector. On the Table \ref{tab}, we show  a sum up  this shift with each corresponding oriented bond and for each oriented bond correspond a type.  
\begin{table}[htp]
\begin{center}
\begin{tabular}{ccccccc}\hline %
Bond on the dual graph && Type && Probability to be open &&Shift\\ \hline
  $\swarrow$ &&  $1$  &&   1  && $(-1,-1)$   \\
	$\longrightarrow$ && $2$   &&  $\alpha$   && $(~~2,~~0)$   \\
  $\nwarrow$ &&  $3$  &&   1  && $(-1,~~1)$   \\ \hline
\end{tabular}
\caption{{\it In this table we describe some elements of the dual graph.} \label{tab}}
\end{center}
\end{table}

We call a {\it nice path}, a path starting at the origin, passing several bonds in the directions of the arrows, loopless and without entries $13$ and $31$. Here as suggested in \cite{Toom:1990}, but do not used there, $123$and $321$ can not occur. Each nice path has a {\it weight}, given by $\alpha^k$ where $k$ is the quantity of $2$ in the path.  For $r\in\{1,2,3\},$ we denote the sum of weights of all the nice paths with $n$ bonds which end in the  vertex $(i,t)$ and have the last bond $r,$ by $S_r(i,t,n).$  Thus, using (\ref{eq:limitesup}) we get
\begin{equation}
\delta_1\Stav^t(x_1=0,\ldots,x_m=0)\le\sum_{n=1}^{\infty}\sum_{r\in\{1,2,3\}}S_r(2m,0,n).
\label{eq:stavsupmeasu}
\end{equation}
Due to the nice path definitions, the numbers $S_r(i,t,n)$ satisfy the initial conditions
\[
S_r(i,t,1)=\left\{\begin{array}{cc}1, &\hbox{ if }i=-1, t=-1 \hbox{ and }r=1, \\ 0, &\hbox{ in all the other cases,}\end{array}\right.
\]

\noindent
and  satisfy the transition equations
\[\left\{\begin{array}{lcl}
S_1(i,t,n+1)&=&S_1(i+1,t+1,n)+S_2(i+1,t+1,n)),\\[.2cm] 
S_2(i,t,n+1)&=&\alpha(S_1(i-2,t,n)+S_2(i-2,t,n)+S_3(i-2,t,n)),\\[.2cm]
S_3(i,t,n+1)&=&S_2(i+1,t-1,n)+S_3(i+1,t-1,n),\\[.2cm]
S_1(i,t,n+2)&=&S_1(i+1,t+1,n+1)+\alpha(S_1(i,t+1,n)+S_2(i,t+1,n)),\\[.2cm] 
S_2(i,t,n+2)&=&\alpha(S_1(i-2,t,n+1)+S_2(i-2,t,n+1)+S_3(i-2,t,n+1)),\\[.2cm]
S_3(i,t,n+2)&=&\alpha(S_2(i,t-1,n)+S_3(i,t-1,n))+S_3(i+1,t-1,n+1). 
\end{array}\right.\]
Let us define sums 
\[
S_r(n)=\sum_{i=-\infty}^{\infty}\sum_{t=-\infty}^{\infty}p^{i}q^{t}S_r(i,t,n), \hbox{ for }r\in\{1,2,3\},
\]
where $p$ and $q$ are non-negative real values.
These quantities  satisfy the initial conditions
\[
S_1(0)=p^{-1}q^{-1},\quad S_2(0)=S_3(0)=0,
\]
\noindent
and recurrence relations for $n\ge 1$,
\[\left\{\begin{array}{lcl}
S_1(n+2)&=&(p^{-2}q^{-2}+\alpha q^{-1})(S_1(n)+S_2(n)),\\[.2cm] 
S_2(n+2)&=&(\alpha pq^{-1}+\alpha^2p^4)S_1(n)+(\alpha pq^{-1}+\alpha^2p^4+\alpha p q)S_2(n)+(\alpha^2p^4+\alpha p q)S_3(n),\\[.2cm]
S_3(n+2)&=&(\alpha q+p^{-2}q^2)(S_2(n)+S_3(n)). 
\end{array}\right.\]
  
 We shall write this system in the  way, $S(n+2)=S(0)M^n,$ where $S(n)=(S_1(n), S_2(n), S_3(n))$ and
\[
M=\left(
\begin{array}{lll}
p^{-2}q^{-2}+\alpha q^{-1}&\alpha pq^{-1}+\alpha^2p^4&0\\
p^{-2}q^{-2}+\alpha q^{-1}&\alpha pq^{-1}+\alpha^2p^4+\alpha p q&\alpha q+p^{-2}q^2\\
0&\alpha^2p^4+\alpha p q&\alpha q+p^{-2}q^2
\end{array}
\right).
\]	
Using (\ref{eq:stavsupmeasu}) and the $S_r(n)$ definition we get
\begin{equation}
\delta_1\Stav^t(x_1=0,\ldots,x_m=0)\le p^{-2m}\sum_{n=1}^{\infty}(S_1(n)+S_2(n)+S_3(n)).
\label{eq:estimations}
\end{equation}

So, if  the right-side of (\ref{eq:estimations}) is convergent, and as $p>1$ we can take $m$ such that this limit will be less than one. 

\subsection{Choice of p and q}
By the Perron-Frobenius theorem, the convergence of (\ref{eq:estimations}) occurs when the maximal eigenvalue of $M$, $\lambda_{pf}$, is less than one. Using a corollary of this theorem (see \cite{Gantmacher}) a necessary and sufficient condition to $\lambda_{pf}\le 1$ is that all the three dominant minors of the matrix $I-M$ are positives. At our case,  verify all these three conditions is too hard. But, we are able to verify the following one among the three
\begin{equation}
1-p^{-2}q^{-2}-\alpha p^{-1}>0 \Rightarrow \alpha<\frac{1-p^{-2}q^{-2}}{p^{-1}}.
\label{eq:pqchoice}
\end{equation}
 Our task is maximize the right-side of (\ref{eq:pqchoice}). As $\alpha$ is less than or equal to one, we get
\[
\frac{1-p^{-2}q^{-2}}{p^{-1}}<1 \Rightarrow -\frac{1}{\sqrt{p(p-1)}}<q<\frac{1}{\sqrt{p(p-1)}}\Rightarrow 1\le q<\frac{1}{\sqrt{p(p-1)}}.
\]
The last implication is a consequence of the fact that $p>1$ and $q\ge 1.$ The inequality 
 $
1\le \frac{1}{\sqrt{p(p-1)}} 
$
is satisfied for $p\in\left(1,\frac{1+\sqrt{5}}{2}\right].$

Now, let us define $f(p,q)=\frac{1-p^{-2}q^{-2}}{p^{-1}}.$
Given $p>1,$ $f(p,q)$ is  increasing as function of $q\ge1$. So,
\[
f(p,q)\le f\left(p,\frac{1}{\sqrt{p(p-1)}}\right)=1.
\]   
Therefore, we finish our task of maximize $f(p,q),$  concluding that  
\begin{equation}p\in\left(1,\frac{1+\sqrt{5}}{2}\right] \hbox{ and }q=\frac{1}{\sqrt{p(p-1)}}. \label{eq:region}\end{equation}
Considering (\ref{eq:region}), we performed some numerical studies for $\lambda_{pf},$ which 
leads us to take $p=\frac{1+\sqrt{5}}{2}.$
 
From the previous statements, we get
\[
\lambda_{pf}= \frac{1}{4}(2\alpha\sqrt{5}+3\alpha^2\sqrt{5}+7\alpha^2+4\alpha+3-\sqrt{5}) +\frac{\sqrt{f(\alpha)}}{4},
\] 
where $f(\alpha)=14+28\alpha^2\sqrt{5}+72\alpha^2+94\alpha^4+4\alpha\sqrt{5}+4\alpha-6\sqrt{5}+172\alpha^3+76\alpha^3\sqrt{5}+42\alpha^4\sqrt{5}.$
For $\alpha\in(0,1),$ the $\lambda_{pf}<1$ when  $\alpha< 0.1142.
$
So, we conclude the proof of Theorem \ref{theor:main}.

\section*{Acknowledgments}
This work has been partially supported by Capes, Coordena\c{c}\~ao de Aperfei\c{c}oamento de Pessoal de N\'ivel Superior - Brasil, under the Program MATH-AMSUD/CAPES (88881.197412/2018-01). P.M.R. has been supported also by FAPESP (2017/10555-0) and CNPq (Grant 304676/2016-0).

\end{document}